\documentclass{article}
\usepackage{amssymb}
\usepackage{graphicx}
\usepackage{amsmath}

\newtheorem{theorem}{Theorem}

\newtheorem{proposition}[theorem]{Proposition}

\newcommand{\SL}{\mathrm{SL}}
\newcommand{\PSL}{\mathrm{PSL}}
\newcommand{\GL}{\mathrm{GL}}
\newcommand{\EL}{\mathrm{EL}}

\newcommand{\Mat}{\mathrm{Mat}}
\newcommand{\Cay}{\mathrm{Cay}}
\newcommand{\Sym}{\mathrm{Sym}}
\newcommand{\Alt}{\mathrm{Alt}}

\newcommand{\mfr}[1]{\mathfrak{#1}}
\newcommand{\la}{\langle}
\newcommand{\ra}{\rangle}
\newcommand{\R}{\mathbb{R}}
\newcommand{\Z}{\mathbb{Z}}
\newcommand{\N}{\mathbb{N}}
\newcommand{\F}{\mathbb{F}}

\newcommand{\ve}{\varepsilon}
\newcommand{\vp}{\varphi}

\begin{document}
\title{Finite Simple Groups as Expanders}

\author{
Martin Kassabov \\
\small{Department of Math} \\
\small{Cornell University}\\
\small{Ithaca, NY 14853}\\
\small{USA.}\\
\and
Alexander Lubotzky%
\thanks{
Supported by NSF and BSF (US--Israel).
This work was done while visiting the IAS in Princeton,
supported by the Ambrose Monell Foundation and
The Ellentuck Fund.
}\\
\small{Institute of Mathematics}\\
\small{Hebrew University}\\
\small{Jerusalem, 91904}\\
\small{Israel.}\\
\and
Nikolay Nikolov\\
\small{New College}\\
\small{Oxford, OX1 3BN}\\
\small{UK.}\\
}
\date{}

\maketitle

\begin{abstract}
We prove that there exist $k\in\N$ and $0<\ve\in \R$ such that every non-abelian finite simple group $G$,
which is not a Suzuki group,
has a set of $k$ generators for which the Cayley  graph $\Cay(G; S)$ is an $\ve$-expander.
\end{abstract}

\section*{Introduction}
\label{sec:intro}

Let $X$ be a finite graph and $0<\ve\in \R$.
The graph $X$  is called an
$\ve$-expander graph if for every subset $A$ of the vertices of $X$
with $|A|\leq \frac{1}{2}|X|, |\partial A|\geq \ve|A|$,
where $\partial A$ denotes the boundary of $A$, i.e., the vertices
of distance $1$ from $A$.

Expander graphs play an important role in computer science and combinatorics
and many efforts has been dedicated to their
constructions (cf.~\cite{expanderbook},the references therein and~\cite{RVW2}).
Many of these constructions are of Cayley graphs and in particular various
families of finite simple groups were shown to be
expanding families.

We will say that a family  $\mfr{G}$ of groups is a family of
expanders (uniformly), if there exists $k\in \N$ and
$0<\ve \in \R$ such that every group $G\in \mfr{G}$,
has a subset $S$ of $k$ generators for which
$\Cay (G; S)$ is an $\ve$-expander.
Until recently all such families were of bounded Lie rank
(cf.~\cite{expanderbook,Ysh}), a fact which has raised the speculation that this
is the only possibility (see~\cite{lubwiess} and~\cite{LR2}).

\medskip

It is easy to see that expander graphs have logarithmic diameter.
In~\cite{BLK} it was shown that every non-abelian finite simple group has a
set $S$ of $7$ generators for which the Cayley
graph $\Cay(G; S)$ has a  diameter bounded by $C\log (|G|)$ for an absolute constant $C$.
It was conjectured there that  one can even make all finite simple groups expanders
in a uniform way although, as observed by
Y.~Luz (see~\cite{lubwiess}), the generators used in~\cite{BLK} do not give rise to expanders.

The main goal of this note is to announce an almost complete proof of this conjecture:

\begin{theorem}
\label{main}
There exist $k\in\N$ and $0<\ve\in\R$ such that every non-abelian finite
simple group $G$, which is not a
Suzuki group, has a set $S$ of $k$ generators for which $\Cay (G; S)$ is an $\ve$-expander.
Moreover the constants $k$ and $\ve$ can be estimated.
\end{theorem}

We believe that the above theorem holds also for the Suzuki groups but our (diverse)
methods do not apply to them.
What makes them exceptional is the fact that they do not contain copies of
$\SL_2(\F_p)$ or $\PSL_2(\F_p)$ like essentially
all the other finite simple groups (see below for more details).

\medskip

The proof of the Theorem is the accumulation of the works~\cite{KNtau,KSL3k,Ksymfull,Nikolov}
and~\cite{Lub2} in the following chronological order:

In~\cite{KNtau} it was proved, extending the work of Shalom~\cite{Ysh},
that for $d\geq 3$ and $k\geq 0$ the group $\SL_d(\Z[x_1, \ldots,x_k])$
has property $(\tau)$, i.e., its finite quotients are expanders.
It was then further shown in~\cite{KSL3k} that if
$R=\Z\la x_1, \ldots, x_k\ra$ is the free non-commutative ring then
suitable finite quotients of $\EL_d(R)$ are also expanders uniformly
(where $\EL_d(R)$ denotes the subgroup of the multiplicative group of $\Mat_d(R)$
generated by the elementary matrices).
This includes, in particular, the groups $\EL_3(\Mat_n(\F_q))\simeq \SL_{3n}(\F_q)$ and thus
$\SL_{3n}(\F_q)$ are uniformly expanders for all
$n$ and for all prime powers $q$.
For every $d\geq 3$, the group $\SL_d(\F_q)$ is a bounded product of copies of
$\SL_{3n}(\F_q)$ for $n=\lfloor d/3 \rfloor$, which implies that the family of groups
$\{\SL_d(\F_q)\,\mid\,d\geq 3, q \mbox{ prime power}\}$
form a family of expanders. 

\medskip

It was then shown in~\cite{Nikolov}  that every finite simple classical group of
Lie type is a bounded product of copies of $\SL_d(\F_q)$ and its central quotients.
This can be combined with the previous results to yield that all classical groups of
Lie type of sufficiently high rank are expanders.

\medskip

The alternating groups $\Alt(n)$  for $n\geq 5$ are also expanders in a uniform way.
This is proved in~\cite{Ksymfull}.
The argument here is more involved since $\Alt(n)$ contains copies of $\SL_d(\F_q)$
but it is not boundedly generated by such groups.
The proof, therefore, goes by decomposing the regular representation of $\Alt(n)$
into two components - on the first $\Alt(n)$ acts ``like''
it is bounded generated by some copies of powers of $\SL_d(\F_q)$
and on the second component the eigenvalues estimates of~\cite{Ro}
are applied (for details, see \S\ref{sec:Alt} below).
This method is influenced by the work of Roichman~\cite{Ro1}.
We are using here the well known representation-theoretic
reformulation of the expanding property - see \S\ref{sec:repth}.

\medskip

On the other hand it is shown in~\cite{Lub2} that the family $\SL_2(\F_q)$,
$q$ a prime power, is also a family of expanders.
This is proved by a combination of Selberg's Theorem (cf.~\cite[\S 4]{expanderbook} )
with the {\it explicit} construction of Ramanujan graphs
as given in~\cite{LSV2}.
In fact, a more careful analysis there can also give the general  result for
the family $\SL_d(\F_q)$ for all $d\geq 2$ and all prime
powers $q$ together, using the theory of Ramanujan complexes (cf.~\cite{LSV1})
and their explicit construction in~\cite{LSV2}.

The case of $\SL_2(\F_q)$ is the crucial new case of~\cite{Lub2}.
It is further shown there, using some model theoretic results of
Hrushovski and Pillay~\cite{HP1}, that for
a fixed $r$, all finite simple groups of Lie type and of rank at most $r$,
with the exception of the Suzuki groups, are bounded
products of copies of $\SL_2(\F_q)$'s.
One therefore deduces that the groups of Lie type and rank $\leq r$ are expanders uniformly.
This exactly complements the results of~\cite{KSL3k} and~\cite{Nikolov} and all together
gives that all finite simple groups of Lie type, with the
possible exception of the Suzuki groups, are expanders uniformly.

\medskip

By the classification of the finite simple groups (CFSG)
every (non-abelian) finite simple group is
either alternating or of Lie type or one of finitely many sporadic groups.
Thus Theorem~\ref{main} follows from the works described  above. \medskip

\medskip

The layout of the current note would not follow the chronological story.
We recall and make in \S\ref{sec:repth} some  observations regarding
the representation theoretic reformulation of the problem.
In \S\ref{sec:SL2} we describe the proof for $\SL_2$, while sections~\ref{sec:SLdRam} and~\ref{sec:SLdT}
give two proofs for the case of $\SL_d, (d\geq 3)$
first via Ramanujan's complexes and second via $\EL_3(\Z\langle x_1, \ldots, x_k\rangle)$.
In section~\ref{sec:Lie} we use the results for $\SL_d(\F_q)$ to construct expanding generating sets in all
simple groups of Lie type (with the exception of the Suzuki groups).
In \S\ref{sec:Alt} we describe the proof for $\Alt(n)$, or equivalently for $\Sym(n)$ -
the symmetric groups, which is a case of great special
interest.

\section{A representation theoretic interpretation}
\label{sec:repth}

As is well known (cf.~\cite{expanderbook} and references therein)
the expanding property of $\Cay(G; S)$ is equivalent to representation
theoretic properties.
The one which will be most convenient for us to use is:

\begin{proposition}
A family $\Cay(G_i; S_i)$ of Cayley graphs are $\ve$-expanders for some $\ve>0$
iff there exist $0<\alpha \in \R$
such that the following holds:

For every $i\in \mathbb{N}$, every unitary  representation
$(V, \vp)$ of $G_i$, every $0\not=v\in V$ and every $0<\delta \in \R$,
if $|\vp(s)v-v |<\delta$ for each $s\in S$, then
$|\vp(g)v-v|<\alpha\delta$ for each $g\in G$.

(i.e., a vector $v$ which is ``$S$-almost invariant'' is also ``$G$-almost invariant''.)
\end{proposition} \medskip

Another interpretation is via the normalized adjacency matrix of the graph $X=\Cay (G; S)$:
Let $k=|S \cup S^{-1}|$ and $\Delta = \frac{1}{k} A$ where $A$ is the adjacency matrix of $X$.
Then all the eigenvalues of $\Delta$ are in $[-1, 1]$.
The second to largest eigenvalue will be denoted $\lambda(X)$ and we have:

\begin{proposition}
\label{eigen}
A family $X_i= \Cay(G_i; S_i)$ is a family of $\ve$-expanders for some
$\ve>0$ if $\lambda(X_i)< 1-\delta$ for some $\delta>0$ and for all
$i$.
\end{proposition}

Note that in Proposition~\ref{eigen} we do not assume that the sets $S_i$ are of bounded size.
If $S_i$ are bounded then the converse is also true. \medskip

Finally, we will make the following observations which will be used several times:
If $S_1$ is a ``bounded product'' of a set $S_2$ i.e.,
every element of $S_1$ is a product of a bounded number of elements of
$S_2$, then an $S_2$-almost invariant vector is also $S_1$-almost invariant.
This implies for example that if $G$ is a bounded product of a bounded number
of subgroups $H_1, H_2, \ldots, H_\ell$ and each
$H_i$ is an expander (w.r.t.\ some set of generators $S_i$) then $G$ is
an expander w.r.t.\ their union  $S=\cup S_i$.

In addition, if $H_i$ are expanders with respect to some generating sets $S_i$,
and each $H_i$ is embedded in a group $G$,
which is an exapnder w.r.t.\ the union $\cup H_i$ then $G$ is
also an expander w.r.t.\ to $\cup S_i$.

\section{The $\SL_2(\F_{p^k})$ case}
\label{sec:SL2}

In this  section we will show,  following~\cite{Lub2},
that $\PSL_2(q)$ can be made into expanders uniformly for every prime power
$q=p^\alpha$.

It is known that: \medskip

\begin{enumerate}
\item
\label{SL2p}
The family of groups $\SL_2(\F_p)$, $p$ prime, with respect to the generators
$A =\left(\begin{array}{cc}1&1\\0&1\end{array}\right)$ and
$B =\left(\begin{array}{cc}0&1\\-1&0\end{array}\right)$
form a family of expanders.
For a proof, based on Selberg's theorem $\lambda_1\geq\frac{3}{16}$, see \cite[\S4]{expanderbook}.

\item
\label{SL2pk}
For a fixed prime $p$, the groups $\SL_2(\F_{p^k})$, $k\in \N$,
have a subset $S_p$ of $p+1$ generators for which the
Cayley groups $X=\Cay(SL_2(p^k); S_p)$ are Ramanujan graphs, i.e.,
$\lambda(X)\leq\frac{2\sqrt p}{p+1}$.
This means, in particular, that they form a family of expanders with a common
expanding factor $\ve>0$, even if we let both $p$
and $k$ vary.
The problem is that the number of generates is unbounded when $p\to\infty$.
\end{enumerate}

The result of \ref{SL2pk} was first proved by Morgenstern~\cite{Mo},
and it relies on the Drinfeld's solution to the characteristic $p$ Ramanujan
conjecture for $\GL_2$.
However, for our purpose we need the explicit construction of Ramanujan graphs
(as special cases of Ramanujan complexes) as given in~\cite{LSV2}.
In the construction there, a symmetric set of $p+1$ generators $S_p'$ for
$\SL_2(\F_{p^k})$ is given as the $p+1$ conjugates of a
fixed element $C$ by a non-split torus $H$ of $\SL_2(\F_p)$, i.e.,
$S_p' =\{h^{-1}C^{\pm 1} h \, \mid\,  h\in H\}$.
Now,
$$
\lambda(\Cay(\SL_2(\F_{p^k}); S_p'))\leq\frac{2\sqrt p}{p+1} < \frac{19}{20}
$$
therefore these graphs are expanders by Proposition~\ref{eigen} uniformly for all $p$
and all $k$.

At the same time $\Cay(SL_2(p);\{A, B\})$ are expanders by~\ref{SL2p}.
Altogether this implies that
$\Cay (\SL_2(\F_{p^k}); \{ A, B, C\})$
are expanders.
Indeed, if $(V, \vp)$ is a representative of
$G=\SL_2(\F_{p^k})$ and $v\in V$ is $S$-almost invariant for the set
$S=\{A, B, C\}$,
then $v$ is $\SL_2(\F_p)$-almost invariant and hence is also $S_p'$-almost invariant,
because $S_p' \subset \SL_2(\F_p) \cdot C \cdot \SL_2(\F_p)$,  and hence $G$-almost
invariant.

\medskip

This shows that $\SL_2(\F_{p^k})$ are uniformly expanders with 3 generators.

\section{$\SL_d(\F_{p^k})$ 
via Ramanujan complexes}
\label{sec:SLdRam}

One can generalize the proof for $\SL_2$ described in \S\ref{sec:SL2} to $\SL_d$ for every $d$.
We sketch the proof below -- for full details see~\cite{Lub2}. \bigskip

Let $F=\F_q$ be the field of order $q=p^k$, for some prime $p$.
Let $E=\F_{q^d}$ be the unique field extension of degree $d$.
The natural map $E^*\to \GL_d(\F_q)$ given by letting $E^*$ act on
$E$ by multiplication, induces an isomorphism of
$$
\{x\in E^*\mid \mbox{Norm}_{E/F}(x)=1\} \mbox{ onto } H\leq \SL_d(\F_q),
$$
where $H$ is a non-split maximal torus of order $\frac{q^d-1}{q-1}$.
In~\cite{LSV2} it was shown that for a suitable choice of
$D\in \SL_d(\F_q)$, the Caley graphs $\Cay(\SL_d(\F_q); S')$ are expanders when
$S'=\{h^{-1} D h\, \mid\, h\in H\}$.
What is really proved there is that the adjacency operator of this Cayley graph is the sum of the
first and the last of the $d-1$ Hecke operators  on the Cayley complex $\SL_d(\F_q)$,
which is a Ramanujan complex, see~\cite{LSV1}.
The eigenvalue bound on these Hecke operators implies that
$$
\lambda(\Cay(\SL_d(\F_q); S'))\leq \frac{2dq^{(d-1)/2}}{(q^d-1)/(q-1)} 
$$
and in particular these are expanders.

\medskip

We mention that the proof here (when $d\geq 3$) does not need the full power  of the Ramanujan bound.
A weaker bound
$$
\lambda(\Cay(\SL_d(\F_q); S'))\leq \frac{1}{\sqrt{q}}+ o(1) \leq \frac{19}{20}
$$
which is sufficient for out purposes,
can be deduced from the infinite dimensional representation theory  of the group
$\SL_d(\F_q((t))\,)$.

\medskip

Assume now that $d$ is even, $d=2m$.
Then the map $E^*\hookrightarrow \GL_d(\F_q)$ described above factors through:
$$
E^*=\F^*_{q^{2m}} \hookrightarrow \GL_2(q^m) \hookrightarrow \GL_{2m}(q) =GL_d(q),
$$
which shows that $H\leq \SL_d(\F_q)$ is contained in a copy of $\SL_2(\F_{q^m})$.

Take now the expanding generating set
$\{A, B, C\}$ of $\SL_2(\F_{q^m})$
from \S~\ref{sec:SL2} with $D$ to
construct a $4$ generator set for
$\SL_d(\F_q)$.

We claim that $\SL_d(\F_q)$ are expanders w.r.t.\ these $4$ generators.
Indeed, again an $S=\{A, B, C, D\}$-almost
invariant vector is $\SL_2(\F_q)$-almost
invariant and then it is $H\cdot D \cdot H$-almost invariant hence also $S'$-almost invariant and thus
$G=\SL_d(\F_q)$-almost invariant as required.

\medskip

This covers the case of even $d$.
For odd $d$, the group $\SL_d(\F_q)$ is a product of 5 copies of
$\SL_{d-1}(\F_q)$ and hence it is also an expander.

\section{$\SL_d(\F_{p^k})$ for $(d\geq 3)$ via $\EL_3(R)$}
\label{sec:SLdT}

A different proof for the case of $\SL_d$ has been given in~\cite{KSL3k}
(prior to the proof in~\cite{Lub2}, which is described in \S\ref{sec:SLdRam}.)
This proof borrows ideas from the work of Shalom~\cite{Ysh}, where he gave a new proof of
Kazhdan's result~\cite{kazhdan} that $\SL_3(\Z)$ has
property $(T)$.
Shalom's proof has two ingredients:
\begin{enumerate}
\item
\label{relT}
The group $\Lambda=\SL_2(\Z) \ltimes \Z^2$ has a relative Kazhdan property (T) relative to
$\Z^2$, i.e., if
$(V, \vp)$ is a unitary representation  space of $\Lambda$ with almost invariant vector,
w.r.t.\ a fixed set of generators, then this vector
is $\Z^2$-almost invariant.

\item
\label{boundedgen}
The group $\Lambda$ has several obvious embeddings into $\SL_3(\Z)$. The copies of
$\Z^2$ in these embeddings (which include the root subgroups) boundedly generate%
\footnote{For our purposes it is sufficient that the quotient groups
$\SL_3(\F_p)$ are uniformly boundedly generated by the elementary matrices.
This result follows easily from the standard row and column reduction algorithm for matrices.}
$\SL_3(\Z)$ by a result of Carter and Keller~\cite{CK}.
\end{enumerate}

Parts \ref{relT} and~\ref{boundedgen} together imply that if $V$ is a unitary representation of
$\Gamma=\SL_3(\Z)$ and $v\in V$ is an almost invariant vector
(w.r.t.\ some generating set) it will be $\Gamma$-almost invariant and therefore the
group $\Gamma$ has property (T).

\medskip

A similar argument to~\ref{relT} works also for
$\Lambda=  \SL_2(R_0) \ltimes R_0^2$, when $R_0$ is the polynomial  ring in $k$ variables over $\Z$.
However, it is not known if the bounded generation of~\ref{boundedgen} holds in this case,
which would imply that
$\SL_3(R_0)$ has property (T).
Still, by analyzing the congruence kernel of $\SL_3(R_0)$, it is deduced in~\cite{KNtau} that
$\Gamma=\SL_3(R_0)$ has property $(\tau)$,
i.e. all of its finite factors are expanders (w.r.t.\ a fixed set of generators of $\Gamma$).

\medskip

A step further is taken in~\cite{KSL3k}: it is shown there that part~\ref{relT}
holds even if one takes the non-commutative free ring
$R=\Z\la x_1, \ldots, x_k\ra$, i.e.,
$\Lambda=R^2\times \EL_2(R)$ has the property (T) relative to $R^2$.
It will be quite surprising if $\Gamma=\EL_3(R)$ is boundedly generated by images of $R^2$ under
different embedding of $\Lambda$ in $\Gamma$, however this holds uniformly for
many quotients $\overline R$ of $R$, probably even for all finite quotients.
For example if
$\overline R=\Mat_n(\F_q)$ or $\overline R=\Mat_n(\F_q)^{\times s}$ for $s\leq q^{n^2}$,
then $\overline R$ is an image of $\Z\la x_1,x_2,x_3\ra$ and
$\EL_3(\overline R)$ is boundedly generated by its elementary matrices
in a way which is independent of $n$, $q$ and $s$.
All this implies that the groups $\EL_3(\Mat_n(\F_q))=\SL_{3n}(\F_q)$ and
$\EL_3(\Mat_n(\F_q)^{\times s})=\SL_{3n}(\F_q)^{\times s}$
are expanders uniformly.

For a general $d\geq 3$, one notes again that $\SL_d(\F_q)$ is a bounded product of a
bounded number of copies of $\SL_{3n}(\F_q)$ for $n=\lfloor d/3 \rfloor$.
We can therefore deduce that $\SL_d(\F_q)$ are expanders for all $d\geq 3$ and every prime power $q$.

One of the advantages of this construction is that it also produces expanders in very
large powers of the group $\SL_{d}(\F_q)$, which are essential for constructing
expanding generating sets in the alternating groups, see \S\ref{sec:Alt}.

\section{Simple groups of Lie type}
\label{sec:Lie}

The following two theorems show how to deduce the (almost) general case of finite simple groups
of Lie type from the case of
$\SL_d(\F_q)$ described in \S\ref{sec:SL2}--\S\ref{sec:SLdT}.

\begin{theorem}[\cite{Nikolov}]
\label{highrank}
There exists a constant $C$ such that every finite simple group $G$ of classical type is a
product of length at most $C$ of at
most $C$ subgroups of $G$ which are quotients of $\SL_d(\F_q)$ (for some $d\geq 2$ and $q$).
\end{theorem}

\begin{theorem}[\cite{Lub2}]
\label{smallrank}
Given $r\in\N$, there exists a constant $C(r)$ with the following property:
Suppose $G$ is
a finite simple group of Lie type of Lie rank at most $r$, which is not a Suzuki group.
Then $G$ is a product of length at most
$C(r)$ of some $C(r)$ subgroups of $G$ which are all quotients of
$\SL_2(\F_q)$ for some $q$'s.
\end{theorem}

As explained in \S\ref{sec:repth}, once a group $G$ is a product of
boundedly  many  expander groups, it is also an expander.
So Theorem~\ref{main} is now proved for all groups of Lie type except for the Suzuki groups.
The reason for the Suzuki groups to be excluded is the fact that they
do not contain copies of $\SL_d(\F_q)$ or $\PSL_d(\F_q)$ for any $d\geq 2$
and prime power $q$.
In fact, the order of a Suzuki group is not divisible by 3 while $|\SL_d(\F_q)|$ is.

The proof of Theorem~\ref{highrank} is by a detailed analysis of the subgroup structure
of simple groups of Lie type: \medskip

The essential point is that arbitrary connected Dynkin diagram becomes
one of type $A_d$ after a few vertices are deleted.
In this way we can find a quasisimple quotient $G_1$ of $\SL_{d+1}(\F_q)$ inside $G$
(in fact it is a Levi factor of a suitable parabolic subgroup).

If $U$ (resp.
$U_1$) is a suitable unipotent subgroup of $G$ (resp $G_1$) then~\cite{Nikolov} proves that
$U$ is a product of at most 14 conjugates of $U_1$.
(using that the positive root system of $G$ parameterizing $U$ is
``close'' to the root system $A_d$ of $U_1$.)
A theorem of M. Liebeck and L. Pyber in~\cite{LiebeckPyber}
now gives that $G$ is a product of at most $13$ conjugates of $U$
and hence a product of at most $13 \times 14=182$ conjugates of $G_1$. \medskip

A very detailed and laborious analysis of this kind
will also lead to a similar proof of Theorem~\ref{smallrank} with explicit bounds for $C(r)$ but
this is not the way this theorem is proved in~\cite{Lub2}.

Instead, the author there appeals to a model theoretic method
developed by Hrushovski and Pillay~\cite{HP1}.
There it is shown that ``definable'' subgroups of $\GL_n(F)$ over
pseudo-algebraically closed field $F$ are very much like Zariski closed
subgroups over algebraically closed field.
In particular, if a definable subgroup $H$ is generated by finitely many
definable subgroups $L_1, \ldots, L_c$ then it is a bounded
product of them.
Now, it follows from the Lang-Weil Theorem that ultra-products of
finite fields are pseudo-algebraically closed (see~\cite{FJ}).
As elementary statements are true in an ultra product
iff they are true in almost all factors -- one can get ``bounded results'' over
finite fields.
This scheme is applied to show that all the finite simple groups
(except the Suzuki groups) contain ``definable'' subgroups isomorphic
to $\SL_2(\F_q)$ or $\PSL_2(\F_q)$ and hence are generated by them  in a bounded
 way (when the rank is bounded).

\section{The Alternating groups}
\label{sec:Alt}
Last but not least is the case of the alternating groups $\Alt(n)$ or
equivalently the symmetric groups $\Sym(n)$.
For its special importance, we restate it:

\begin{theorem}[\cite{Ksymfull}]
There exist $k\in \N$ and $0<\ve\in\R$ such that for every
$n\in\N$ the symmetric group $\Sym(n)$ has an explicit set of
generators $S_n$ of size at most $k$ such that $\Cay(\Sym(n); S_n)$
is an $\ve$-expander.
\end{theorem}

The main idea of the proof in~\cite{Ksymfull} is as follows:
Let first consider the alternating groups $\Alt(n)$. Assume first that
$n=d^6$ and $d=2^{3k}-1$ for some $k$.
Based on ideas and results from~\cite{KSL3k} (see \S\ref{sec:SLdT} above)
it is shown first
that the Cayley graphs of the groups $\Delta_k = \SL_{3k}(\F_2)^{d^5}$ with respect to some
generating set $F_k$ of size at most $20$ are expanders for all $k$ and
$d=2^{3k}-1$. \medskip

Now, thinking of the set $\{1, \ldots, n\}$ as the points in a $6$-dimensional cube
of size $d$, and remembering that $\SL_{3k}(\F_2)$ acts
transitively on a set of size $d$ (via its defining linear action on $\F_2^{3k}$),
we can construct $6$ different embeddings $\pi_i$ of $\Delta_k$ into
$\Alt(n)$, where the image under $\pi_i$ of each copy of $\SL_{3k}(\F_2)$ in $\Delta_k$
acts on the points on a line parallel to the $i$-th coordinate axis in the cube.
Denote $S_n=\cup_i\, \pi_i(F_k)$ and $E=\cup_{i}\,\pi_i(\Delta_k)$.
We will show that the Cayley graphs of $\Alt(n)$ with respect to
$S_n$ form a family of expanders. Using that $F_k$ is an expanding generating set
in $\Delta_k$ it suffices to show that the existence of $E$-almost invariant
vector $v$ in any unitary representation $V$ of $\Alt(n)$ implies the existence of
an invariant vector.

We decompose $V$ as a sum of two representations $V_1 \oplus V_2$, where
$V_1$ is the sum of all irreducible representations of $\Alt(n)$ in $V$ corresponding to
partitions $\lambda\vdash n$, where the first row $\lambda_1$ is small (less than $n - d^{5/4}$)
and $V_2$ contains all others. A similar decomposition is used by Roichman in~\cite{Ro1}
to show that the Cayley graphs of $\Alt(n)$ with respect to conjugacy classes have
some expanding properties. We will use two different arguments to show that the projection $v_1$ of $v$ in
$V_1$ is small and that projection $v_2$ in $V_2$ is close to an invariant vector:

\begin{enumerate}
\item
Using the definition of the set $E$,
it is shown that a bounded power of $E$ acts transitively on ``nearly all''
ordered $\ell$-tuples of points in the cube for some $\ell$ of size approximately $d^5/(3\log d)$.
Also a bounded power of $E$ contains a permutation which acts as $\ell$-cycle.
Thus the vector $v$, and therefore $v_1$, is almost invariant by ``nearly all'' elements in the conjugacy
class $C_\ell$ of $\ell$-cycles in $\Alt(n)$. Here ``nearly all'' means a subset of proportion
tending to $1$, as $k$ and $n$ tend to infinity.
This implies that $v_1$ is almost invariant under the action of the operator
$$
L=\frac{1}{|C_\ell|} \sum_{s\in C_\ell} s.
$$

The operator $L$ acts as a multiplication by $\chi_\lambda(C_\ell) / \chi_\lambda(id)$ on the
irreducible representation $V_\lambda$ corresponding to the partition $\lambda$, where
$\chi_\lambda$ is the character of $V_\lambda$.
At this point one can appeal to methods and results of Roichman~\cite{Ro},
who studied normalized character values of the symmetric groups.
Roichman's results gives that $|\chi_\lambda(C_\ell) / \chi_\lambda(id)| \ll 1$, for
any  $\lambda\vdash n$, provided that the first row $\lambda_1$ is small. Therefore
$|| L v_1 || \ll ||v_1||$,
which together with the almost invariance of $v_1$ under $L$,
implies that the vector $v_1$ is short.

\item
Using that all irreducible representations in $V_2$ corresponds to partitions with
$\lambda_1 \geq n- d^{5/4}$,
one can view the linear span $W$ of the orbit of $v_2$ in $V_2$ as part of the induced
representation to $\Alt(n)$ of the trivial representation of $\Alt(n- d^{5/4})$.
This induced representation has a basis $\mfr{B}$, whose elements corresponds to
the ordered $d^{5/4}$-tuples of points in the cube.
The size of the basis $\mfr{B}$ is significantly smaller that the size of $E$ and it
can be shown that the random walk on $\mfr{B}$ defined by $E$ mixes in just several steps.
This together with the almost invariance of $v_2$ under $E$,
implies that $v_2$ is close to an invariant vector.
\end{enumerate}

This finishes the sketch of the proof for the case $\Alt(n)$ for $n=(2^{3k}-1)^6$ for some $k$. The case of
general $n$ follows from the observation that $\Alt(n)$ can be written as
a product of a bounded number of copies of $\Alt(n_k)$ for $n_k=(2^{3k}-1)^6$
embedded in $\Alt(n)$. Finally, any expanding generating set of
$\Alt(n)$ gives an expanding generating set in $\Sym(n)$ by adjoining one odd permutation.

\bibliographystyle{amsplain}
\bibliography{FSG}

\end{document}